\documentclass[12pt]{article}
\usepackage{amsmath, amsthm, amssymb}
\usepackage{hyperref}

\newtheorem{thm}{Theorem}
\newtheorem{lem}[thm]{Lemma}

\numberwithin{equation}{section}
\begin{document}
\title{Convergence of Ricci flow on $\mathbb{R}^2$ to flat space}
\author{James Isenberg\thanks{Partially supported by the NSF under Grant 
 PHY-652903 }
\\ University of Oregon
\\ \small{isenberg@uoregon.edu}  \and
Mohammad Javaheri 
\\ Trinity College
\\ \small{Mohammad.Javaheri@trincoll.edu}}

\date{15 April, 2009}
%\author{James Isenberg \\ \emph{\small{isenberg@uoregon.edu}}\\ Mohammad Javaheri\\ \emph{\small{Mohammad.Javaheri@trincoll.edu}}}
\maketitle

\begin{abstract}
We prove that, starting at an initial metric $g(0)=e^{2u_0}(dx^2+dy^2)$ on $\mathbb{R}^2$ with bounded scalar curvature and bounded $u_0$, the Ricci flow $\partial_t g(t)=-R_{g(t)}g(t)$ converges to a flat metric  on $\mathbb{R}^2$. 
\end{abstract}

\section{Introduction}
In two dimensions, the Ricci flow equation ${{\partial} \over {\partial t}} g(t)=-2Ric_{g(t)} $ reduces to 
\begin{equation}
\label{firsteq}
{{\partial} \over {\partial t}} g(t)=-R_{g(t)}g(t),
\end{equation}
where $Ric_{g(t)}$ denotes the Ricci curvature of $g(t)$, and  $R_{g(t)}$ denotes its  scalar curvature. We are interested in the long-term behavior on $\mathbb{R}^2$ of conformally flat solutions to the Ricci flow equation. For these, the metrics take the form $g(x,t)=e^{2u(x,t)}g_E$, where $g_E$ is the standard Euclidean metric on $ \mathbb{R}^2$. A straightforward  calculation shows that if we set $v(x,t):=e^{2u(x,t)}$, then the initial value problem for the flow equation \eqref{firsteq} takes the form 
\begin{equation} \label{secondeq}
{\partial \over {\partial t}}v(x,t) =\Delta \ln v(x,t)~,~v(x,0)=e^{2u_0(x)}~.
\end{equation}
which is the ``fast diffusion" initial value problem on  $\mathbb{R}^2$

The long-term existence of solutions of \eqref{firsteq} or \eqref{secondeq} has been studied in \cite{das}, where it is shown that the solutions exist for all $t\geq 0$ if and only if
\begin{equation} \label{intcond}
\int_{\mathbb{R}^2}e^{2u_0}dxdy =\infty~.
\end{equation}
The following theorem by L. Wu \cite{lfw} addresses the long-term behavior of the solutions of \eqref{firsteq}. We first need a few definitions. Let $g$ be a complete metric on $\mathbb{R}^2$. The \emph{circumference at infinity} is defined as
\begin{equation}C_\infty(g)=\sup_K \inf_{D} \{L(\partial D)|\forall~ \mbox{compact set}~K\subset \mathbb{R}^2, \forall ~\mbox{open set}~ D \supset K \}~,\end{equation}
and the \emph{aperture} of $g$ is defined as
\begin{equation}A(g)={1 \over {2\pi}} \lim_{r \rightarrow \infty} {{L(\partial B_r)} \over r}~.\end{equation}
Here $B_r$ denotes the geodesic ball (or disc) of radius $r$ and $L(\partial B_r)$ is the length of the boundary of $B_r$. Also the Ricci flow is said to have \emph{modified subsequence convergence}, if there exists a 1-parameter family of diffeomorphisms $\{\phi_t\}_{t \geq 0}$ such that for any sequence $t_i \rightarrow \infty$, there exists a subsequence (denoted again by $t_i$) such that the sequence $\phi_{t_i}^*g_{t_i}$ converges uniformly on every compact set as $i \rightarrow \infty$. In terms of these definitions, Wu's theorem can be stated as follows: 
\begin{thm} \label{wu} \cite{lfw}
Let $g(t)=e^{2u(t)}g_{E}$ be a solution to \eqref{firsteq} such that $g(0)=e^{2u_0}g_E$ is a complete metric with bounded curvature and $|\nabla u_0|$ is uniformly bounded on $\mathbb{R}$. Then the Ricci flow has modified subsequence convergence as $t\rightarrow \infty$. Moreover if the curvature is positive at time zero, then the limiting metric is a cigar soliton\footnote{The cigar soliton on $\mathbb{R}^2$ has the metric given by \eqref{cigar} below.} if $C_\infty(g(0))<\infty$, or a flat metric if $A(g(0))>0$.     \end{thm}

Our main result is to prove a similar theorem which replaces the condition that the curvature of $g(0)$ be positive with the condition that $u$ be bounded throughout $\mathbb{R}^2$. Specifically, we shall prove the following:

\begin{thm}\label{limit0}
Suppose $g_0=e^{2u_0}g_E$ has bounded curvature and $u_0$ is a bounded
smooth function on $\mathbb{R}^2$. Then the Ricci flow $\partial_t
g=-Rg$ exists for all $t\geq 0$ and has modified subsequence convergence to the flat metric in the $C^k$-
topology of metrics on compact domains in $\mathbb{R}^2$ for each
$k$.
\end{thm}

The boundedness condition on $u_0(x)$ guarantees that the initial metric $g_0$ is complete, and also guarantees that the condition \eqref{intcond} holds. The boundedness condition on the curvature, combined with standard elliptic estimates implies that $|\nabla u(x,0)|$ is bounded. Note that  if $u_0$ is bounded, then the aperture of $g(0)$ is positive; however, a priori it is possible that $R_{g(0)}$ is negative on some region in $\mathbb{R}^2$. We will show through Lemmas \ref{bounds} and \ref{boundR} that $R_{g(t)} \rightarrow 0$ uniformly as $t\rightarrow \infty$. This together with similar asymptotic bounds which we show hold for the covariant derivatives of the curvature will prove Theorem \ref{limit0}. 

Theorem \ref{limit0} is sharp in the sense that there are example solutions of Ricci flow for which
$u$ is not bounded initially and there is no subsequence for which
$g(t)$ converges to a flat metric. In fact the cigar soliton on $\mathbb{R}^2$ given by
\begin{equation}\label{cigar}
g_\Sigma(x,t)=\left ({1 \over { |x|^2+e^{4 t}}} \right )g_E
\end{equation}
provides such a example solution. In this example $u_0=-{1 \over 2}\ln(1+|x|^2)\leq 0$, but $u_0$ has no lower bound. The results in \cite{ses-das} show that if $g(t)=v(x,t)g_E$ is a solution of the Ricci flow such that the scalar curvature stays bounded and the width of $g(t)$ is finite for each $t$, then $g(t)$ is a gradient soliton of the form \eqref{cigar} up to translations and scalings. Here the width of a metric $g$ on $\mathbb{R}^2$ is defined as 
\begin{equation}w(g)=\inf_F \sup_c L\{F=c\}~,~F:\mathbb{R}^2 \rightarrow [0,\infty)~,\end{equation}
and $L\{F=c\}$ is the length of the level curve $F=c$ in the metric $g$. 

The work of Hsu \cite{hsu1,hsu2} shows that there are alternative sets of conditions which one can place on solutions of \eqref{secondeq} which also guarantee that the solutions approach cigar soliton geometries. To state Hsu's theorem\footnote{We thank the referee for telling us about Hsu's results.},  let us set
$$\phi_{\beta, k}={2 \over {\beta(|x|^2+k)}}~.$$ Then one has the following:

\begin{thm} \label {hsu} \cite{hsu1,hsu2}
Suppose $v(x,t)$ is a solution of \eqref{secondeq} so that 
$$\mbox{i)}~ \liminf_{|x| \rightarrow \infty} {{\ln v(x,t)} \over {\ln |x|}} \geq -2 ~\mbox{uniformly in}~[t_1,t_2]~~\forall t_2>t_1>0~,\hspace{400pt}$$
$$\mbox{ii)}~{\partial \over {\partial t}}v(x,t) \leq {{v(x,t)} \over t}~,\hspace{400pt}$$
$$\mbox{iii)}~\phi_{\beta, k_1} \leq v(x,0) \leq \phi_{\beta, k_2}~,\hspace{400pt}$$
for some $\beta>0$ and $k_1>k_2>0$. Then the solution $v(x,t)$ of \eqref{secondeq} will converge uniformly on $\mathbb{R}^2$ and also in $L^1(\mathbb{R}^2)$ to the function $\phi_{\beta, k}$ as $t\rightarrow \infty$ for some $k>0$ which is uniquely determined by the initial function $v(x,0)$. 
\end{thm}

In the spirit of Theorems \ref{limit0} and \ref{hsu}, it is rather tempting to make the following conjecture by replacing the conditions (i) and (ii) of Theorem \ref{hsu} by a geometric condition involving only the initial geometry:
\\
\\
\textbf{Conjecture}. \emph{Let $g(t)=v(x,t) g_E$ be such that $g(0)$ has bounded scalar curvature and there are $\beta, k_1,k_2>0$ such that $\phi_{\beta, k_1} \leq v(x,0) \leq \phi_{\beta, k_2}$. Then the flow exists for all times and converges to $\phi_{\beta, k}g_E$ for some $k>0$, uniformly on compact sets in the $C^m$-topology of metrics for each $m$. }

\section{Maximum Principles}

There are a number of different versions of the Maximum Principle, corresponding to a variety of different parabolic or elliptic partial differential equations on a number of different domains, with various added hypotheses applied to the solutions. We use two versions here. The first, which we label MP1, states that if a $w(x,t)$ is a solution  to the heat equation ${\partial \over {\partial t}} w = \Delta_{g(t)} w $  on $\mathbb{R}^2 \times (0,\infty)$, then $|w(x,t)| \leq |w(x,0)|$ for all $t \in (0,\infty)$. This is standard, and needs no proof here. Since the second version (stated below and labeled MP2) is less standard, so we do prove it here. First, we need the following lemma, which is proved in \cite{lfw}:
\begin{lem}\label{lem13} \cite[Lemma 1.3]{lfw}
There exists a non-negative time-independent function
$\eta:\mathbb{R}^2 \times [0,T]\rightarrow \mathbb{R}$ such that
$\eta(0,t)=0$, $\Delta \eta \leq 1$ and
$\eta(x,t)=\eta(|x|)\rightarrow \infty$ as $|x|\rightarrow \infty$, for all $t\in [0,T]$.
\end{lem}

We now state and prove 
\begin{thm} \textbf{\emph{(MP2)}} \label{mp}
Let $g(t)$ be a smooth one-parameter family of metrics on
$\mathbb{R}^2$ for $t\in [0,T]$. Suppose a function
$u:\mathbb{R}^2 \times [0,T]\rightarrow \mathbb{R}$ satisfies
\begin{equation}{\partial \over {\partial t}} u \leq \Delta_{g(t)} u + hu~,\end{equation}
where $h$ is a smooth function. If $u\leq 0$ at $t=0$ and
both $u$ and $h$ are uniformly bounded on $\mathbb{R}^2 \times [0,T]$, then $u(x,t)\leq
0,~\forall (x,t) \in \mathbb{R}^2 \times [0,T]$.
\end{thm}

\begin{proof}
Choose $C>0$ such that $|h(x,t)|< C$ for all $(x,t)\in
\mathbb{R}^2\times[0,T]$. Define $\bar u:=e^{-Ct}u$. Then, by
setting $\bar h:=h-C<0$, we have
\begin{equation}{\partial \over {\partial t}} \bar u \leq \Delta_{g(t)} \bar u +\bar h\bar u~.\end{equation}
Let $\eta$ be the function obtained in Lemma \ref{lem13}. For each $\delta>0$, let
\begin{equation}u_\delta := \bar u - \delta \eta -\delta t~.\end{equation}
Then we verify that 
\begin{equation}\label{udelta}
    {\partial \over {\partial t}} u_\delta \leq \Delta u_\delta +\delta(\Delta
\eta-1)+\bar hu_\delta +\bar h\delta(\eta+t)~.
\end{equation}
Since $\bar u$ is bounded and $\eta(|x|)\rightarrow \infty$ as
$|x|\rightarrow \infty$, the maximum of $u_\delta$ is attained at
some $(x_0,t_0)\in \mathbb{R}^2 \times [0,T]$. We show that this
maximum is non-positive. Otherwise if $u_\delta(x_0,t_0)>0$, then
$t_0>0$. On the other hand, the inequality \eqref{udelta} implies that
\begin{equation}0\leq {\partial \over {\partial t}} u_\delta(x_0,t_0)\leq \bar h\delta(\eta+t_0)~,\end{equation}
which is a contradiction, since $\bar h<0$. Hence $u_\delta(x,t)
\leq 0$ for all $(x,t) \in \mathbb{R}^2 \times [0,T]$. By letting
$\delta \rightarrow 0$ we conclude that $\bar u\leq 0$ and so
$u\leq 0$ for all $(x,t)$, which in turn implies that $u\leq 0$ for all $(x,t) \in \mathbb{R}^2 \times [0,T]$.
\end{proof}

\section{Convergence to the flat metric}

We prove our main result Theorem \ref{limit0} here, preceded by a sequence of lemmas which establish  the estimates we need to complete the proof. 

In accordance with the hypotheses of Theorem \ref{limit0}, we assume here  that $g_0=e^{2u_0}g_E$ has bounded scalar
curvature $|R_0|<k_0$ and that $u_0(x)=u(x,0)$ is bounded. As noted above, it  follows from standard elliptic gradient estimates (see for example Theorem 3.9 \cite{GT}) that $|\nabla u_0|$ is bounded on $\mathbb{R}^2$. Let $g(t)=e^{2u(x,t)}g_E$ be
the Ricci flow starting at $g(0)=g_0$. The long-term existence of the flow follows from \cite{das} or Theorem \ref{wu}. Replacing the quantity $u(x,t)$ for the moment by 
\begin{equation}f(x,t):=-2u(x,t)\end{equation}
we have the initial value problem
\begin{equation}\label{propf}
   {\partial \over {\partial t}} f=\Delta_{g(t)} f=R_{g(t)}~,~f(x,0)=f_0(x)~.
\end{equation}

Applying Theorem 2.4 from reference \cite{lfw} to this flow, we obtain\footnote{Note that the hypotheses for Theorem 2.4 in \cite{lfw} do \emph{not} include a positivity condition on the curvature} a uniform bound on $R(x,t)$ as well as a uniform bound on $|\nabla f(x,t)|$. It then follows from \eqref{propf} and MP1 that $f(x,t)$ is uniformly bounded for all $(x,t) \in \mathbb{R}^2 \times [0,\infty)$. We can in fact improve the bound on the scalar curvature as follows:

\begin{lem}\label{bounds}
Suppose $g(0)$ has bounded scalar curvature and $u_0$ is bounded. Let $g(t)$ be the solution of the Ricci Flow $\partial_t g=-Rg$ on $\mathbb{R}^2$. Then there exists a constant $C>0$ depending only on $g(0)$ such that:
\begin{equation}-{C /({1+Ct})} \leq R(x,t) \leq C~,~\forall (x,t) \in \mathbb{R}^2 \times [0,\infty)~.\end{equation}
\end{lem}
\begin{proof}Choose $k_0$ such that $|R(x,0)|\leq k_0$. 

Let
\begin{equation}\theta(t):={{-k_0} \over {1+k_0t}}~\end{equation}
and define $S(x,t):=R(x,t)-\theta(t)$. Then one easily calculates that $S(x,t)$ satisfies the PDE
\begin{equation}{\partial \over {\partial t}}S=\Delta S +(R+\theta)S~.\end{equation}
Since $S(x,t)$ and $R(x,t)+\theta(t)$ are uniformly bounded and since by definition 
$S(x,0)\geq 0$, Theorem \ref{mp} (MP2) applied to $-S(x,t)$ implies that
$S(x,t)\geq 0$ for all $x,t$. It follows that $R(x,t) \geq
-k_0/(1+k_0t)$ and the proof of Lemma \ref{bounds} is
completed.
\end{proof}

Our next task is to obtain a better control on
$|\nabla f|$ as $t\rightarrow \infty$.

\begin{lem}\label{boundnabla}
Let $g_0=e^{2u_0}g_E$ be a complete metric on $\mathbb{R}^2$.
Suppose $u_0$ and the scalar curvature of $g_0$ are bounded. Then
there exists a constant $C$ that depends only on $g_0$ such that
\begin{equation} \label{nablaf}\sup_{x\in \mathbb{R}^2} |\nabla f(x,t)|^2 \leq {C \over {1+t}}~.\end{equation}
\end{lem}
\begin{proof}
Set
\begin{equation}F:=t|\nabla f|^2+f^2~.\end{equation}
The evolution equation \eqref{propf} and the evolution equation for $|\nabla f|^2$ under the Ricci flow
\begin{equation}\label{evolnablaf}
    {\partial \over {\partial t}} |\nabla f|^2 = \Delta |\nabla f|^2-2|D^2f|^2\leq \Delta |\nabla f|^2~
\end{equation}
imply that ${\partial_t}F
\leq \Delta F$. Since $F(x,t)$ is uniformly bounded on $\mathbb{R}^2 \times [0,T]$ for all $T>0$, the Maximum Principle MP1 implies that $F(x,t)$ is uniformly bounded on $\mathbb{R}^2 \times [0,\infty)$, and 
recalling that $f(x,t)$ is also bounded, we obtain \eqref{nablaf}.
\end{proof}

We next let 
\begin{equation}\label{defH}
    H:=R(x,t)+|\nabla f(x,t)|^2~, 
\end{equation}
and readily verify that $H(x,t)$ satisfies
\begin{equation}\label{eveq}
    {\partial \over {\partial t}} H=\Delta H - 2\left |D^2 f -{1 \over 2}\Delta f \cdot
g \right |^2~.
\end{equation}
We can then show

\begin{lem}\label{boundR}
Under the same hypothesis as in Lemma \ref{boundnabla}, we have
\begin{equation}\sup_{x\in \mathbb{R}^2}H(x,t)\leq {C \over {1+t}}~,\end{equation}
for some $C>0$ depending only on $g(0)$.
\end{lem}
\noindent \emph{Sketch of the proof}. We verify this estimate for $H(x,t)$  following the ideas used in proving \cite[Prop. 5.30]{cd}.
Let $G:=t(H+|\nabla f|^2)$. Using simple estimates, one can show that $\partial_t G \leq
\Delta G$ for $t$ large enough. Since $G(x,t)$ is uniformly bounded on $\mathbb{R}^2 \times [0,T]$ for each $T>0$, MP1 implies that $G(x,t)$ is uniformly bounded for all $(x,t) \in \mathbb{R}^2 \times [0,\infty)$.  \hfill $\square$

It follows from the definition of $H(x,t)$ and from this estimate that 
 $R(x,t) \leq H(x,t)  \leq C/(1+t)$. Comparing this with the
result from Prop. \ref{bounds}, we conclude that $R(x,t)\rightarrow 0$
as $t\rightarrow \infty$.

We now obtain estimates for derivatives of the curvature.  First, we show
\begin{lem}
Under the same hypotheses as in Lemma \ref{boundnabla}, we have
\begin{equation} \label{derivR}\sup_{x\in \mathbb{R}}|\nabla R(x,t)|^2 \leq {C \over
{(1+t)^3}}~,~\forall t \geq 0~.\end{equation}
\end{lem}

\begin{proof} Let $J:=t^4|\nabla R|^2+\lambda t^3R^2$, where the constant 
$\lambda$ is to be chosen later. Since $\partial_t |\nabla R|^2
\leq \Delta |\nabla R|^2+4R|\nabla R|^2$ and $\partial_t R=\Delta
R+R^2$, we have
\begin{eqnarray}\nonumber
% \nonumber to remove numbering (before each equation)
   \partial_t J &=& t^4 \partial_t |\nabla R|^2+4t^3|\nabla
R|^2+\lambda t^3\partial_t R^2+3\lambda t^2R^2\\ \nonumber
   &\leq &  \Delta J +4t^3(tR+1-\lambda/2)|\nabla
   R|^2+\lambda(2t^3R^3+3t^2R^2)~.
\end{eqnarray}
As a consequence of Lemma \ref{boundR} one can choose $\lambda$ so that $tR+1\leq
\lambda/2$ ; hence
\begin{equation}\partial_t J \leq \Delta J+C~,\end{equation}
for some constant $C$. On the other hand, it follows from Shi's derivative estimates \cite{shi} that we have
\begin{equation}\label{shis}
|\nabla^m R(x,t)| \leq {K_{m,k_0} \over {\sqrt{t^m}}}~,~\forall m \geq 0,~\forall t\in (0,T_{k_0}]~,
\end{equation}
where $K_{m,k_0}$ and $T_{k_0}$ are constants depending only on $m,k_0$. This clearly implies that $|\nabla R| \leq C$ for some constant $C$ and $t$ large enough. In particular $J(x,t)$ is uniformly bounded for all $(x,t) \in \mathbb{R}^2 \times [0,T]$ for any $T>0$. It again follows from the Maximum Principle (MP2) that $J(x,t)$ is uniformly bounded for all $(x,t) \in \mathbb{R}^2 \times [0,\infty)$, and \eqref{derivR} follows. .
\end{proof}

\begin{lem}
 Under the same hypothesis as in Lemma \ref{boundnabla}, for each $k$ there exists
a constant $C_k$ that depends only on $g_0$ such that
\begin{equation}\label{higherbounds}
    \sup_{x\in \mathbb{R}^2}|\nabla^k R(x,t)|^2 \leq {{C_k} \over
{(1+t)^{k+2}}}~.
\end{equation}
\end{lem}
\begin{proof}
The argument which proceeds from the estimates established for the curvature $R(x,t)$ and for its first derivative in the above lemmas to the estimates \eqref{higherbounds} for all orders of derivatives of the curvature closely follows the pattern of the proof of the same result for geometries on closed 2 manifolds, as discussed in  \cite[Prop. 5.33]{cd}. The only real difference is in the use of the maximum principle. Here, we use MP2, along with Shi's derivative estimates \eqref{shis}
 
\end{proof}

We are now ready to complete the proof of Theorem \ref{limit0}.
\\
\emph{Proof of Theorem \ref{limit0}}. Since $f(x,t)$ is uniformly bounded by a constant $K>0$ on $\mathbb{R}^2 \times [0,\infty)$, for any fixed $x_0\in \mathbb{R}^2$ and for any sequence $t_i\rightarrow \infty$, there exists a subsequence, denoted by $t_i$ again, such that $\beta=\lim_{i \rightarrow \infty} f(x_0,t_{i})$ exists. Allowing $x\in
\mathbb{R}^2$ to vary, we can establish (for all $t\geq0$)
\begin{equation}|f(x,t)-f(x_0,t)|\leq d_t(x,x_0)\sup_{x\in \mathbb{R}^2}|\nabla f(x,t)|
\leq e^{K/2}d(x,x_0)\sqrt{{C \over {1+t}}}~,\end{equation} where
$d_t(x,x_0)$ is the distance between $x$ and $x_0$ in
$g(t)$ and $d(x,x_0)$ is the usual Euclidean
distance. It follows that $f(x,t_{i})$ is also convergent to
$\beta$ as $i\rightarrow \infty$. In other words
\begin{equation}g(t_{i})=e^{2u(x,t_{i})}g_E\rightarrow e^{-\beta}g_E~,~i\rightarrow \infty~.\end{equation}
 The $C^k$-convergence of $g(t_{i})$ to this flat limit then follows from
\eqref{higherbounds}. \hfill $\square$

\end{document}